\documentclass[12pt]{amsart}
\usepackage{amsmath,amssymb
}
\theoremstyle{plain}
\newtheorem{theorem}{Theorem}
\newtheorem{lemma}[theorem]{Lemma}

\newtheorem{proposition}[theorem]{Proposition}

\newtheorem{aclaim}{Claim} 

\theoremstyle{remark}

\newtheorem*{ack}{Acknowledgement}

\theoremstyle{definition}

\def\endpfclaim#1{
\renewcommand{\qedsymbol}{\openbox${}_{\;\mathrm{Claim#1}}$}
\end{proof}
\renewcommand{\qedsymbol}{\openbox}
}

\def\bplus{b_+^2}
\def\phu{\varphi}
\def \R {\mathbf{R}}

\def\EE{\mathcal{E}}

\def\LL{\mathcal{L}}

\def\MM{\mathcal{M}}

 \def\MM{\mathcal{M}}

\def\SC{\ifmmode{\text{SPIN}^c}\else{$\text{SPIN}^c$}\fi}

\DeclareMathOperator{\spinc}{Spin^c}

\def\sw{Seiberg--Witten}
\def\SW{\ifmmode{\text{SW}}\else{$\text{SW}$}\fi}
\def\SWW{\ifmmode{{\text{SW}}_{tot}}\else{${\text{SW}}_{tot}$}\fi}

\begin{document}
\title[A vanishing theorem]{The Seiberg--Witten
invariants of manifolds with wells of negative curvature}
\author[Daniel Ruberman]{Daniel Ruberman}
\address{Department of Mathematics, MS 050\newline\indent Brandeis
University \newline\indent Waltham, MA}
\email{\rm{ruberman@brandeis.edu}}
\date\today \maketitle
\section{Introduction}
A $4$-manifold with $\bplus >1$ and a nonvanishing Seiberg--Witten invariant
cannot admit a metric of positive scalar curvature.  This remarkable fact is
proved~\cite{witten:monopole} using the Weitzenb\"ock--Lichnerowicz formula for the square of the
$\spinc$ Dirac operator, combined with the `curvature' part of the \sw \
equations. Thus, in dimension $4$, there is a strong generalization of
Lichnerowicz's vanishing
theorem~\cite{lichnerowicz:spinors,lawson-michelson,rosenberg-stolz:psc} for
the index of the Dirac operator of a spin manifold with a metric of positive
scalar curvature.

In the recent years, the method of semigroup domination
\cite{err:curvature,rosenberg:semigroup,rosenberg:vanishing} has led to a
different sort of generalization of Lichnerowicz's  theorem and other theorems
in which a positive curvature hypothesis  leads to a topological vanishing
theorem.  Essentially, the  hypothesis of positive curvature may be weakened to
permit negative curvature on a `small' set.  (The precise notion of `small'
depends on what kind of curvature is being discussed--see the statement of
Theorem~\ref{swvanish} for the version we are using.)   In this note, we use  semigroup domination to show
that a
$4$-manifold  with positive scalar curvature away from a set of small volume must have
vanishing  Seiberg--Witten invariants.  Moreover, the same vanishing holds for the \sw \
invariant of any finite covering space.

We sketch very briefly the definition of the \sw\ invariants, and refer to~\cite{morgan:swbook,nicolaescu:swbook,moore:swbook} for more details.   Recall that a $\spinc$ structure
$\sigma$ on a smooth Riemannian $4$-manifold $X$  determines a pair of spinor
bundles $W^\pm \to X$ which are Hermitian bundles over $X$ of rank $2$.  A unitary connection $A$ on $\LL = \det(W^+)$
determines the Dirac operator
$$
D_A^+ : \Gamma(W^+) \xrightarrow{\nabla_A} \Gamma(T^*X  \otimes W^+)
\xrightarrow{\rho}\Gamma(W^-)
$$
where $\nabla_A$ is the induced connection on $W^+$ and $\rho$
denotes Clifford multiplication.  The \sw\ equations, for a
connection $A$ and spinor $\phu$ are
\begin{equation}\label{sw}
\left\{\begin{aligned}
D_A^+\phu &= 0\\
\rho(F_A^+) +i\mu^+ &= i\tau(\phu,\phu)
\end{aligned}
\right\}
\end{equation}
Here $\tau(\phu,\phu) = (\phu\otimes\phu^*)_0$ denotes the traceless
part of the endomorphism $\phu\otimes\phu^*$ of $W^+$, and $\mu$ is a
real $2$-form.   The solution space to equations~\eqref{sw},
modulo gauge equivalence, gives the \sw \ moduli space $\MM(X,\sigma)$.  It is compact, and for a generic choice of $(g,\mu)$ is an oriented smooth manifold of dimension $$
\dim(\sigma) = \mathrm{ind}_\R (D_A^+) -\frac{1}{2}(e(X) + \mathrm{sign}(X))
$$
If $d=\dim(\sigma)$ is negative or odd, then the \sw\ invariant is defined to be $0$.  Otherwise, the \sw\ invariant is defined as $\SW_X(\sigma) = \langle c_1(\EE)^{d/2},[\MM(X,\sigma)]\rangle$, where $\EE \to \MM$ is a naturally defined complex line bundle.  If $\bplus(X) > 1$, then this count is independent of the metric, while if $\bplus(X) =1$, there is a mild dependence on the metric.  Briefly, in this case, the space of (metrics, $2$-forms) is divided by codimension-one `walls' into connected components called chambers, on which the the invariant is constant. 

To state the main result, let $\MM(n,K,D,V)$ be the class of Riemannian
$n$-manifolds with Ricci curvature $\ge K$, diameter $\le D$ and
volume  $\ge V$.
\begin{theorem}\label{swvanish}
Let $(X^4,g) \in \MM(4,K,D,V)$ have $\bplus> 1$, and let $s_0 > 0$.  There exists
an
$\epsilon> 0$, depending only on $K,D,V$, such that if the scalar
curvature $s(X,g)$ is bounded below by $s_0$ except on a set of
volume less than $\epsilon$, then $\SW_X(\sigma)= 0$ for all $\spinc$
structures on $X$.  Moreover, for any finite cover $\pi: \tilde{X}\rightarrow X$
and
$\spinc$ structure $\tilde{\sigma} = \pi^*\sigma$ on $\tilde{X}$, we
have $\SW_{\tilde{X}}(\tilde\sigma)= 0$.  If $\bplus(X) = 1$, then these
statements hold for the \sw \ invariant associated to the chamber containing
$g$.
\end{theorem}
One can deduce a statement about $\spinc$ structures
$\tilde\sigma$ which do not necessarily pull back from a $\spinc$ structure on
$X$.  This requires an additional hypothesis, and will be discussed
following the proof of the statement.     Note that the same
$\epsilon$ works for finite covers of $X$ with arbitrary degree.   In principle,
a vanishing result for $\SW_{\tilde X}$ could be deduced directly from the
first statement of the theorem.  However, in this approach, the constant
$\epsilon$ (which depends~\cite[Proposition 1.2]{err:curvature} on the volume
and diameter) would decrease rapidly as the degree of the cover increases.
\begin{ack}
The author was partially supported by NSF Grant 9971802.  The
work on this paper was done while the author was visiting the
University of
Paris, Orsay, under the partial support of the CNRS.  Thanks to Steve Rosenberg for introducing me to these sorts of questions and for some helpful comments.
\end{ack}

\section{A vanishing theorem}
The key ingredient in the proof of Theorem~\ref{swvanish} is a generalization
of Kato's inequality, referred to as semigroup domination
\cite{donnelly-li:semigroup,h-s-u:domination,rosenberg:semigroup}. For a
Riemannian metric on $X$, let $\Delta$ denote the Laplacian on functions, and
$s$ the scalar curvature function.   Likewise,
$\tilde\Delta$ and $\tilde s$ will represent the Laplacian and scalar
curvature on a finite covering space $\tilde X$.  Let $L$ be the self-adjoint operator $\nabla_A^*\nabla_A +s/4$, where $\nabla_A$ is the connection on the bundle of spinors, and s is the scalar curvature of the metric.  The operator $L$ is that part of the Weitzenb\"ock decomposition of $D_A^*D_A$ which does not involve the curvature of the bundle $W^+$.   As above, for a $\spinc$ connection $\tilde A$ on an arbitrary $\spinc$ structure on a finite cover of $X$, we will denote by $\tilde L$ the operator $\tilde\nabla_{\tilde A}^*\tilde\nabla_{\tilde A} +\tilde{s}/4$.

The principle of semigroup domination, in its strong form~\cite{rosenberg:vanishing}, yields the following pointwise inequality.
\begin{lemma}\label{semigroup}  For any unitary connection $A$ on $\LL$, and any
$\psi \in L^2(W^+)$
\begin{equation}\label{kato}
\left| e^{-tL}(\psi)(x) \right| \leq  2 e^{-t(\Delta + s/4)} |\psi|(x)
\end{equation}
\end{lemma}
The smoothing property of the operators $e^{-tL}$ and $e^{-t(\Delta + s/4)}$ implies that the quantities being compared in~\eqref{kato} are continuous functions of $x$, so that the inequality makes sense pointwise.  It is also important to note that the right side of this inequality is independent of $A$.

The main step in Theorem~\ref{swvanish} is the following vanishing result.  From now on, the phrase `for sufficiently small $\epsilon$' will be used as a shorthand for the hypotheses
of Theorem~\ref{swvanish}.
\begin{proposition}\label{reducible}
For sufficiently small $\epsilon$, any solution to the SW equations
$$D_A^+\phu = 0,
F_A^+ = \tau(\phu,\phu)$$
is reducible (ie has
$\phu\equiv 0$).
\end{proposition}

\begin{proof}[Proof of Proposition~\protect{\ref{reducible}}]
Proposition 2.4
of~\cite{rosenberg:semigroup}, and the discussion following it shows
\begin{aclaim}\label{positive}
For sufficiently small $\epsilon$, the operator $\Delta +s/4$ is positive, as is $\tilde\Delta
+\tilde s/4$ for every covering space of $X$.
\end{aclaim}

The second step is analogous to~\cite[Theorem 2.2]{rosenberg:semigroup}
\begin{aclaim}\label{eigen}
For sufficiently small $\epsilon$, the operator $L$ has only positive eigenvalues.
\end{aclaim}
\begin{proof}[Proof of Claim~\ref{eigen}]
Suppose that $L\psi = \lambda \psi$ with $\lambda \leq 0$ and $\psi \neq 0$.
By elliptic regularity for $L$, we may assume that $\psi$ is continuous.  In particular, there is a point $x_0$ with $\psi(x_0) \neq 0$.  According to Claim~\ref{positive}, the right hand side of~\eqref{kato}, evaluated at $x_0$,
goes to 0 as
$t
\to
\infty$.  However, the left hand side is
$e^{-t\lambda}|\psi(x_0)|$, which goes to infinity if $\lambda <0$. If $\lambda =0$, then the left hand side is the non-zero constant $|\psi(x_0)|$, so we still get a contradiction.
\endpfclaim{~\ref{eigen}}
\begin{aclaim}\label{vareigen}
For $\epsilon$ sufficiently small, there is no non--zero solution to $L\phu = -f(x)\phu$ where f is
a continuous function with $f(x) \geq 0$ for all X.
\end{aclaim}
\begin{proof}
Choose an orthonormal basis $\{\psi_i\}$ of eigenspinors for L with eigenvalues
$\lambda_i$, which are all positive according to Claim~\ref{eigen}.  Then
$\phu = \sum a_i \psi_i$, where by a standard argument the convergence is in $L^2_2$.   It follows that
$ L\phu = \sum a_i \lambda_i\psi_i$, so that
$$
<L\phu,\phu>_{L^2} = \sum \lambda_i a_i^2 > 0
$$
But
$$
<L\phu,\phu>_{L^2} = - \left(\int f(x) <\phu,\phu>\right)\ \leq 0
$$
which is a contradiction.
\endpfclaim{~\ref{vareigen}}

To prove Proposition~\ref{reducible}, assume that $(A,\phu)$ is a solution to the \sw\ equations.
By regularity for solutions of the \sw\ equations, $\phu$ is smooth, and so is $|\phu|^2$.  Using the Weitzenb\"ock formula, and substituting for the curvature term,
$$ 0= (D^+_A)^*D_A^+\phu = (\nabla_A^*\nabla_A +s/4)\phu +
\frac{1}{2}\tau(\phu,\phu)\phu = L\phu + \frac{1}{4}|\phu|^2\phu
$$
So $L\phu =- \frac{1}{4}|\phu|^2\phu$ and we conclude that $\phu \equiv 0
$ by Claim~\ref{vareigen}.
\end{proof}

\begin{proof}[Proof of Theorem~\protect{\ref{swvanish}}] If $\bplus > 1$, then
(\cite{donaldson-kronheimer},\cite[\S 6.3]{morgan:swbook}) the set of `generic'
metrics for which the  \sw\ equations admit no reducible solution is open and
dense in the space of all metrics (using the $C^\infty$ topology).

For sufficiently small $\delta > 0$, Proposition~\ref{reducible} provides an
$\epsilon >0$ such that if $(X,g') \in \MM(K-\delta,D+\delta,V-\delta)$, and
the scalar curvature $s(X,g') \ge s_0$ except on a set of volume less than
$\epsilon + \delta$, then any solution to the \sw \ equations on $X$ is
reducible.  Now if $(X,g) \in \MM(4,K,D,V)$, and
$$
vol_{g}\{x \in X | s_g(x) < s_0\} \leq \epsilon,
$$
approximate $g$ by a generic metric $g' \in \MM(K-\delta,D+\delta,V-\delta)$
with
$$
vol_{g'}\{x \in X | s_{g'}(x) < s_0\} \leq \epsilon + \delta.
$$
Now any solution to the \sw \ equations is reducible, but  $g'$ is chosen so
that there  are no reducible solutions either.   Since we may compute it with
respect to any metric, the \sw\ invariant must vanish.

To prove the vanishing statement for the \sw\ invariant on a covering space
$\tilde X$, we use the observation of~\cite{rosenberg:semigroup,err:curvature}
that the curvature assumptions on $X$ imply that the operator $\tilde\Delta +
\tilde s/4$ is positive.   Thus the proof of Proposition~\ref{reducible}
applies to show that any solution to the \sw\ equations on $\tilde X$ must be
reducible.   This argument does not use the fact that $\tilde\sigma$ is pulled
back from $X$, but merely that the Laplacian on functions is pulled back from
$X$.   This, in turn only requires that the metric on $\tilde X$ be the
pullback metric.

Now we make use of a simple principle about reducible solutions to the \sw \
equations.
\begin{aclaim}\label{cover}
Suppose that $\tilde \sigma = \pi^*\sigma$.  Then for a generic metric on $X$,
there are no reducible solutions to the \sw\ equations on $\tilde X$.
\end{aclaim}
\begin{proof}[Proof of Claim~\ref{cover}]
Denote by $\tilde \LL= \pi^*\LL$ the determinant bundle of the $\spinc$
structure $\tilde\sigma$.   For a generic metric $g$ on $X$, there is a
$g$--self-dual form $\alpha \in H^2_+(X)$ with $\alpha \cup c_1(\LL) \neq 0$.
This is the content of the generic metrics theorem for $X$.   But, pulling back
to $\tilde X$, this means that $\pi^*\alpha \cup c_1(\tilde \LL) \neq 0$, so
$c_1(\tilde \LL) $ can't be represented by a $\pi^*g$--anti-self-dual form.   Hence
there are no reducible solutions on $\tilde X$.
\endpfclaim{~\ref{cover}}
The rest of the proof of the vanishing theorem works exactly as above.
\end{proof}

We can still give a vanishing result for the \sw\ invariants of a covering
space, even for $\spinc$ structures which don't pull back from $X$.  The result
requires an additional mild topological hypothesis; there are perhaps some
variations of this method giving similar results.
\begin{theorem}\label{swcover}
Let $(X^4,g) \in \MM(4,K,D,V)$ have $\bplus> 1$,  and let $s_0 > 0$.  Assume in
addition that
$$
\frac{1}{2}(e(X) + \mathrm{sign}(X)) = \bplus(X) -b^1(X) +1 > 0
$$
Then there exists
an $\epsilon> 0$, depending only on $K,D,V$, such that if the scalar
curvature $s(X,g)$ is bounded below by $s_0$ except on a set of
volume less than $\epsilon$, then $\SW_{\tilde{X}}(\tilde\sigma)= 0$ for all
$\spinc$ structures $\tilde\sigma$ on $\tilde X$.
\end{theorem}
\begin{proof}
As in the proof given above, for sufficiently small $\epsilon$, and any metric
on $X$ satisfying the hypotheses, any solution to the \sw\ equations on $\tilde
X$ will be reducible.  Let $(\tilde A,0)$ be such a solution, and let
$D^+_{\tilde A}$ be the corresponding Dirac operator.  The Weitzenb\"ock formula
says that
$$
((D^+_{\tilde A})^* D_{\tilde A}) = \tilde L + \rho(F_{\tilde A}^+)
$$
But since $\tilde A$ is reducible, $F_{\tilde A}^+ = 0$, and so $ (D^+_{\tilde A})^*D^+_{\tilde A}$ is the operator $\tilde L$.  But the argument above says that
$\tilde L$ is positive (for sufficiently small $\epsilon$, of course.)  It
follows that $\ker(D^+_{\tilde A}) = 0$, and therefore that the index of
$D^+_{\tilde A}$ is less than or equal to $0$.

By definition, the \sw\ invariant $\SW_{\tilde X}(\tilde \sigma)$ is non-zero only when the dimension of
the \sw\ moduli space $\MM(\tilde X,\tilde \sigma)$ is $\ge 0$.  This dimension is given by the formula
$$
\begin{aligned}
\dim(\tilde\sigma) &= \mathrm{ind}(D^+_{\tilde A}) -
\frac{1}{2}(e(\tilde X) + \mathrm{sign}(\tilde X))\\
&= \mathrm{ind}(D^+_{\tilde A})  -n\cdot\frac{1}{2}(e(X) + \mathrm{sign}(X))
\end{aligned}
$$
Since the first term on the right hand side is non-positive, and we have
assumed that the second term is negative, this gives a contradiction.
\end{proof}
Note that the hypotheses of Theorem~\ref{swcover} are satisfied for $X$ simply
connected.  So the argument presented above gives a simpler way
to prove Theorem~\ref{swvanish} in that case.

The reader may wonder if there are examples of manifolds with vanishing \sw\ invariants to which we can apply the part of Theorem~\ref{swvanish} dealing with covering spaces.  A class of such manifolds was described by Shuguang Wang~\cite{wang:vanishing}; his examples are complex surfaces $\tilde X$ admitting a free anti-holomorphic involution $\tau$.  He shows that the quotient $X = \tilde X /\tau$ has vanishing \sw \ invariants; however Theorems~\ref{swvanish} and~\ref{swcover} still apply.  

\section{Higher dimensions}
There are more elaborate versions of Lichnerowicz' index-theoretic obstruction to the existence of a metric of positive scalar curvature on a spin manifold--see~\cite{rosenberg-stolz:psc} for a recent overview.  It was noted by Hitchin~\cite{hitchin:spinors} that the Dirac operator on a spin manifold $M^n$ defines an element $\alpha(M) \in KO_n$, which vanishes if $M$ admits a metric of positive scalar curvature.   For $n\ge 5$, it was shown by Stolz~\cite{stolz:psc} that a simply-connected spin manifold $M^n$ for which $\alpha(M) =0$ admits a metric of positive scalar curvature.  The construction of positive scalar curvature metrics uses the surgery method of Schoen--Yau~\cite{schoen-yau:psc} and Gromov--Lawson~\cite{gromov-lawson:psc}.

Using these results, we show that a metric on a manifold with `mostly' positive curvature can sometimes be traded in for a metric having everywhere positive curvature.  
\begin{theorem}\label{surgery} Suppose that $n\ge 5$.  Let $(X^n,g) \in \MM(n,K,D,V)$ be simply connected, and let $s_0 > 0$.  There exists an $\epsilon> 0$, depending only on $n,K,D,V$, such that if the scalar curvature $s(X,g)$ is bounded below by $s_0$ except on a set of
volume less than $\epsilon$, then $X$ admits a metric of strictly positive curvature.
\end{theorem}
\begin{proof}
Following~\cite{rosenberg:semigroup}, choose $\epsilon$ sufficiently small so that the operator $\Delta +s/4$ is positive.  Then semigroup domination implies that the kernel and cokernel of the Dirac operator must vanish.  This implies that the invariant $\alpha(X)$ is $0$.  By the result of Stolz, it follows that $X$ must in fact admit a metric of positive scalar curvature. 
\end{proof}
It would be more satisfactory to have a direct method for modifying a metric to eliminate a small region of negative curvature rather than to have to appeal to the theorem of Stolz.  Such an argument would have a better chance of applying in dimension $4$.  Apparently there are results along these lines for Ricci curvature, but none for scalar curvature.  

The argument above should extend to non-simply-connected manifolds whose groups satisfy the `Gromov--Lawson--Rosenberg' conjecture.  This states that the existence of a positive scalar curvature metric on $M$ is equivalent to the vanishing of an index-theoretic invariant $\alpha(M,f) \in KO_n(C^*\pi)$.  Here $f:M \to B\pi$ classifies the universal cover of the fundamental group $\pi$ of $M$ and $KO_n(C^*\pi)$ denotes the K-theory of a $C^*$--algebra associated to $\pi$.  The proof of the required vanishing result would require that semigroup domination applies to the Dirac operator with coefficients in the Hilbert space $C^*(\pi)$.

\providecommand{\bysame}{\leavevmode\hbox to3em{\hrulefill}\thinspace}


\begin{thebibliography}{10}

\bibitem{donaldson-kronheimer}
S.K. Donaldson and P.B. Kronheimer, \emph{{The Geometry of Four-Manifolds}},
  Clarendon Press, Oxford, 1990.

\bibitem{donnelly-li:semigroup}
Harold Donnelly and Peter Li, \emph{Lower bounds for the eigenvalues of
  {R}iemannian manifolds}, Michigan Math. J. \textbf{29} (1982), no.~2,
  149--161.

\bibitem{err:curvature}
{Elworthy D. and Rosenberg, S., with appendix by Ruberman, D.}, \emph{Manifolds
  with wells of negative {R}icci curvature}, Inventiones Math. \textbf{103}
  (1991), 491--496.

\bibitem{gromov-lawson:psc}
Mikhael Gromov and H.~Blaine Lawson, Jr., \emph{The classification of simply
  connected manifolds of positive scalar curvature}, Ann. of Math. (2)
  \textbf{111} (1980), no.~3, 423--434.

\bibitem{h-s-u:domination}
H.~Hess, R.~Schrader, and D.~A. Uhlenbrock, \emph{Domination of semigroups and
  generalization of {K}ato's inequality}, Duke Math. J. \textbf{44} (1977),
  no.~4, 893--904.

\bibitem{hitchin:spinors}
N.J. Hitchin, \emph{Harmonic spinors}, Adv.\ in Math. \textbf{14} (1974),
  1--55.

\bibitem{lawson-michelson}
H.~Blaine Lawson, Jr. and Marie-Louise Michelsohn, \emph{Spin geometry},
  Princeton University Press, Princeton, NJ, 1989.

\bibitem{lichnerowicz:spinors}
Andr{\'e} Lichnerowicz, \emph{Spineurs harmoniques}, C. R. Acad. Sci. Paris
  \textbf{257} (1963), 7--9.

\bibitem{moore:swbook}
John~Douglas Moore, \emph{Lectures on {S}eiberg-{W}itten invariants}, second
  ed., Springer-Verlag, Berlin, 2001.

\bibitem{morgan:swbook}
John~W. Morgan, \emph{The {S}eiberg-{W}itten equations and applications to the
  topology of smooth four-manifolds}, Mathematical Notes, vol.~44, Princeton
  University Press, Princeton, NJ, 1996.

\bibitem{nicolaescu:swbook}
Liviu~I. Nicolaescu, \emph{Notes on {S}eiberg-{W}itten theory}, American
  Mathematical Society, Providence, RI, 2000.

\bibitem{rosenberg-stolz:psc}
Jonathan Rosenberg and Stephan Stolz, \emph{Metrics of positive scalar
  curvature and connections with surgery}, Surveys on surgery theory, Vol. 2,
  Princeton Univ. Press, Princeton, NJ, 2001, pp.~353--386.

\bibitem{rosenberg:vanishing}
Steven Rosenberg, \emph{Semigroup domination and vanishing theorems}, Geometry
  of random motion (Ithaca, N.Y., 1987), Amer. Math. Soc., Providence, RI,
  1988, pp.~287--302.

\bibitem{rosenberg:semigroup}
\bysame, \emph{Applications of semigroup domination}, Diffusion processes and
  related problems in analysis, Vol.\ I (Evanston, IL, 1989), Birkh\"auser
  Boston, Boston, MA, 1990, pp.~285--291.

\bibitem{schoen-yau:psc}
R.~Schoen and S.~T. Yau, \emph{On the structure of manifolds with positive
  scalar curvature}, Manuscripta Math. \textbf{28} (1979), no.~1-3, 159--183.

\bibitem{stolz:psc}
S.~Stolz, {\em Simply connected manifolds of positive scalar curvature}, Ann.
  of Math. (2), {\bf 136} (1992), 511--540.

\bibitem{wang:vanishing}
Shuguang Wang, \emph{A vanishing theorem for {S}eiberg-{W}itten invariants},
  Math. Res. Lett. \textbf{2} (1995), no.~3, 305--310.

\bibitem{witten:monopole}
E.~Witten, \emph{Monopoles and four-manifolds}, Math.\ Res.\ Lett. \textbf{1}
  (1994), no.~6, 769--796.

\end{thebibliography}
\end{document}